\documentclass[12pt,reqno]{amsart}
\usepackage{setspace}
\usepackage{amsmath}
\usepackage{amsfonts}
\usepackage{amssymb}
\usepackage{amscd}
\usepackage{graphics}
\usepackage{verbatim}
\def\a{\alpha}
\def\b{\beta}

\def\bR{{\mathbb R}}

\def\bC{{\mathbb C}}

\def\b1{{\rm id}}

\newfont{\goth}{eufm10 scaled \magstep1}

\newfont{\mcal}{eusm10 scaled \magstep1}

\newtheorem{Th}{Theorem}
\newtheorem{remar}[Th]{Remark}
\newtheorem{Prop}[Th]{Proposition}
\newtheorem{Cor}[Th]{Corollary}
\newtheorem{Lem}[Th]{Lemma}
\newtheorem{Def}[Th]{Definition}
\newtheorem{Ex}[Th]{Example }

\def\bt{\begin{Th}}
\def\et{\end{Th}}
\def\bp{\begin{Prop}}
\def\ep{\end{Prop}}
\def\bc{\begin{Cor}}
\def\ec{\end{Cor}}
\def\bl{\begin{Lem}}
\def\el{\end{Lem}}
\def\bd{\begin{Def}}
\def\ed{\end{Def}}
\def\bex{\begin{Ex}}
\def\eex{\end{Ex}}
\def\br{\begin{remar}}
\def\er{\end{remar}}

\def\be{\begin{equation}}
\def\ee{\end{equation}}
\def\ben{\begin{enumerate}}
 \def\een{\end{enumerate}}
\def\ba{\begin{array}{rlll}}
\def\ea{\end{array}}
\def\bea{\begin{eqnarray}}
\def\eea{\end{eqnarray}}
\def\bean{\begin{eqnarray*}}
\def\eean{\end{eqnarray*}}

\def\diag{\mathrm{diag\;}}

\def\C1{cohomogeneity one  }

\newcommand{\K}{K\"ahler}

\newcommand{\ov}[1]{\overline{#1}}

\newcommand{\Boc}{\operatorname{Boc}}
\newcommand{\Cut}{\operatorname{Cut}}

\newcommand{\Sp}{\operatorname{Sp}}

\begin{document}

\title[Bochner coordinates on flag manifolds]{Bochner coordinates on flag manifolds}

\author{Andrea Loi}
\address{(Andrea Loi) Dipartimento di Matematica \\
         Universit\`a di Cagliari (Italy)}
         \email{loi@unica.it}

\author{Roberto Mossa}
\address{(Roberto Mossa) Dipartimento di Matematica \\
         Universit\`a di Florianopolis (Brasil)}
         \email{roberto.mossa@gmail.com}

\author{Fabio Zuddas}
\address{(Fabio Zuddas) Dipartimento di Matematica e Informatica \\
          Via delle Scienze 206 \\
         Udine (Italy)}
\email{fabio.zuddas@uniud.it}

\thanks{
The first  and third authors were  supported by Prin 2015 -- Real and Complex Manifolds; Geometry, Topology and Harmonic Analysis -- Italy and also by INdAM. GNSAGA - Gruppo Nazionale per le Strutture Algebriche, Geometriche e le loro Applicazioni.}
\subjclass[2000]{53D05;  53C55;  53D05; 53D45} 
\keywords{Symplectic maps; \K\ manifolds; Gromov width; Gromov-Witten invariants}

\begin{abstract}
We find necessary and sufficient conditions under which the complex coordinates on a flag manifold of a classical group described in \cite{A-P} are Bochner coordinates.
\end{abstract}
 
\maketitle

\tableofcontents

\section{Introduction and statement of the main result}
Let us recall that, given a complex manifold $M$ endowed with a real analytic \footnote{Throughout the paper all \K\ manifolds will be
assumed to be real analytic.}
 \K \ metric $g$, Calabi introduced,
in a neighborhood of a point 
$p\in M$,
a very special
\K\ potential
$D_p$ for the metric
$g$, which 
he christened 
{\em diastasis}.
More precisely, recall that a \K\ potential
is an
analytic function 
$\Phi$ 
defined in a neighborhoood
of a point $p$
such that 
$\omega =\frac{i}{2}\partial \bar\partial\Phi$,
where $\omega$
is the \K\ form
associated
to $g$.
In a complex coordinate system
$(z)$ around $p$
$$g_{\alpha\beta}=
2g(\frac{\partial}{\partial z_{\alpha}},
\frac{\partial}{\partial \bar z_{\beta}})
=\frac{{\partial}^2\Phi}
{\partial z_{\alpha}\partial\bar z_{\beta}}.$$
 
A \K\ potential is not unique:
it is defined up to the sum with
the real part of a holomorphic function.
By duplicating the variables $z$ and $\bar z$
a potential $\Phi$ can be complex analytically
continued to a function 
$\tilde\Phi$ defined in a neighborhood
$U$ of the diagonal containing
$(p, \bar p)\in M\times\bar M$
(here $\bar M$ denotes the manifold
conjugated of $M$).
The {\em diastasis function} is the 
\K\ potential $D_p$
around $p$ defined by
$$D_p(q)=\tilde\Phi (q, \bar q)+
\tilde\Phi (p, \bar p)-\tilde\Phi (p, \bar q)-
\tilde\Phi (q, \bar p).$$
Among all the potentials the diastasis
is characterized by the fact that 
in every coordinates system 
$(z)$ centered in $p$
$$D_p(z, \bar z)=\sum _{|j|, |k|\geq 0}
a_{jk}z^j\bar z^k,$$
with 
$a_{j 0}=a_{0 j}=0$
for all multi-indices
$j$.

One can always find local 
(complex) coordinates
in a neighborhood of 
$p$
such that
$$D_p(z, \bar z)=|z|^2+
\sum _{|j|, |k|\geq 2}
b_{jk}z^j\bar z^k,$$
where 
$D_p$ is the diastasis
relative to $p$.
These coordinates,
uniquely defined up
to a unitary transformation,
are called 
{\em the Bochner   or normal coordinates}
with respect to the point $p$ (cfr. \cite{note, boc, ca, ruan,tian4}).
Calabi's  diastasis function has shown to be  the right tool in the study of  the Riemannian geometry of a  \K\ manifold, since it  keeps track of its underlying  complex and symplectic structure, in contrast with the geodesic distance which has the disadvantage to be preserved by maps between Riemannian manifolds only in the totally geodesic case, while the diastasis is preserved by any \K\ map (see \cite{ca} for a proof). 
For example, in  \cite{diastexp} the first two authors of the present paper have defined  the diastatic exponential,
by twisting the classical riemannian one with Calabi's diastasis and obtaining applications   to the symplectic geometry of symmetric spaces, while in  \cite{diastbal, diastent} the second author has introduced and analyzed the concept of diastatic entropy for hyperbolic manifolds
by  extending that of volume entropy using the diastasis function instead of the geodesic distance.

Given a point $p$ of a \K\ manifold $M$, it is then natural to analyze the relationship between the maximal domain of definition  of  Bochner coordinates $\Boc_p$, the  maximal domain   of definition  of the diastasis $D_p$, say $V_p$,  and the cut locus $\Cut_p$ of $p$ (obviously 
$\Boc_p\subset V_p$ by the very definition of Bochner coordinates). 
When the \K\ manifold involved is a Hermitian symmetric space of compact type, H. Tasaki \cite{tas} has shown that  
$\Boc_p=V_p=M\setminus \Cut_p$ (see also \cite{diastherm} for other properties of diastasis function of symmetric spaces and \cite{berezin, remrkhom} for the  case of bounded homogeneous domains).

Therefore if one considers the more general class of  flag manifolds the following question  arises.

\vskip 0.2cm

\noindent
{\bf Question 1:} {\em Given $p$ a point on a flag manifold is it true that $\Boc_p=V_p$?}

\vskip 0.2cm

Regarding the second  equality, we believe the validity of the following:

\vskip 0.2cm

\noindent
{\bf Conjecture:} {\em Assume that $V_p=\Omega\setminus \Cut_p$, for some point $p$ of a flag manifold $\Omega$. Then $\Omega$
is a Hermitian symmetric space of compact type.}
 
\vskip 0.3cm

In Theorem \ref{thm main} below, which represents the main result of the  paper,  we provide a partial answer to Question 1 
(cfr. Remark \ref{remarteor}).
Before stating Theorem \ref{thm main}
let us recall  some basic facts about flag manifolds.

A flag manifold $\Omega = G/K$ of a semisimple compact Lie group $G$ is an orbit of the adjoint action of $G$ on its Lie algebra $\mathfrak{g}$. 
The interest for these manifolds can be understood by recalling (see, for example, \cite{B}) that each compact homogeneous \K\ manifold
$M$ is the \K\ product of a flat complex torus and a simply-connected compact
homogeneous \K\ manifold, and admits a \K-Einstein structure if and only if
is a torus or is simply-connected. In the simply-connected case, $M$ is isomorphic, as
a homogeneous complex manifold, to an orbit of the adjoint action of its connected
group of isometries $G$ (which, being compact and with no center, is semisimple).

A flag manifold endowed with an invariant complex structure can be combinatorially described by a so-called {\it painted Dynkin diagram}, that is the Dynkin diagram of $G$ where some nodes have been painted in black. The Dynkin diagrams of the classical simple groups $SU(d)$, $Sp(d)$, $SO(n)$ are the following

\begin{align*}
&&& \underset{\substack{}}{\circ} - \underset{\substack{}}{\circ} - \dotsb - \underset{\substack{}}{\circ} - \underset{\substack{}}{\circ} &&
 (G= SU(d))
\end{align*}

\begin{align*}
&&& \underset{\substack{}}{\circ} - \underset{\substack{}}{\circ} - \dotsb - \underset{\substack{}}{\circ} \Leftarrow \underset{\substack{}}{\circ} && (G= Sp(d))
\end{align*}

\begin{align*}
&&& \underset{\substack{}}{\circ} - \underset{\substack{}}{\circ} - \dotsb - \underset{\substack{}}{\overset{\overset{\textstyle\circ_{}}{\textstyle\vert}}{\circ}} \,-\, \underset{\substack{}}{\circ} && (G= SO(2d))
\end{align*}

\begin{align*}
&&& \underset{\substack{}}{\circ} - \underset{\substack{}}{\circ} - \dotsb - \underset{\substack{}}{\circ} \Rightarrow \underset{\substack{}}{\circ} && (G= SO(2d+1))
\end{align*}

The number of black nodes in a painted diagram equals the second Betti number $b_2(\Omega)$ of the corresponding manifold $\Omega$.

Moreover, if $\alpha_1, \dots, \alpha_p$ are the black nodes, then an invariant \K\ metric on $G/K$ is determined by a choice of positive real numbers $c_1, \dots, c_p$ associated to $\alpha_1, \dots, \alpha_p$ (see Section 2 below for more details).

When $G$ is one of the classical semisimple groups, then D. V. Alekseevsky and A. M. Perelomov \cite{A-P} describe explicit complex coordinates on a dense open  subset $U\subset \Omega=G/K$ around a point $0\in\Omega$ such that $U$ is biholomorphic to 
$\bC^N$, where $N$ is the complex dimension of $\Omega$. Throughout the paper we are going to call  these coordinates {\it Alekseevsky-Perelomov coordinates}. They also describe  an explicit \K\ potential  $D_0:U\rightarrow\bR$ (see (\ref{potential}) below)  for any \K\ form $\omega$ on $\Omega$ (Theorem \ref{potential} below).
In \cite{LMZ}, the authors of the present paper have proven that the potential  $D_0$ is indeed  Calabi's diastasis function centered in $0$ and moreover $U$ is the maximal domain of definition of $D_0$. 
Hence it makes sense to analyze when  Alekseevsky-Perelomov coordinates are  Bochner coordinates.
Our main result in this direction is the following:
\begin{Th}\label{thm main} Let $\Omega = G/K$ be an irreducible flag manifold of classical type with $b_2(\Omega)=p$, endowed with the invariant \K\ form $\omega$ determined by coefficients $c_1, \dots, c_p > 0$ associated to the black nodes $\alpha_1, \dots, \alpha_p$ of its painted diagram.
Then, the Alekseevsky-Perelomov coordinates are Bochner, up to rescaling
\footnote{Notice that if $D_p(z, \bar z)= \sum_{\alpha} c_{\alpha} |z_{\alpha}|^2+
\sum _{|j|, |k|\geq 2}
b_{jk}z^j\bar z^k$, then clearly it is enough to apply a change of coordinates of the kind $z_{\alpha} \mapsto \lambda_{\alpha} z_{\alpha}$ ($ \lambda_{\alpha} \neq 0$) in order to have Bochner coordinates: in this case, we say that the $z_{\alpha}$'s are {\it Bochner up to rescaling}.}, only in the following cases:
\begin{enumerate} 
\item[(i)] $p = 1$, for every $G$ and every $\omega$;
\item[(ii)] $p = 2$, $G = SU(d)$ and $c_1 = c_2$;
\item[(iii)] $p = 2$, $G=SO(2d)$, the painted diagram of $\Omega$  is

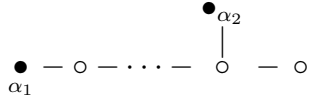
\begin{align*}
&&& \underset{\substack{\alpha_1}}{\bullet} - \underset{\substack{}}{\circ} - \dotsb - \underset{\substack{}}{\overset{\overset{\textstyle\bullet_{\alpha_2}}{\textstyle\vert}}{\circ}} \,-\, \underset{\substack{}}{\circ} &&
\end{align*}
and $c_1 = 2 c_2$.
\end{enumerate} 
\end{Th}
(We recall that if $b_2(\Omega)=1$ then the invariant \K\ forms on $\Omega$ are all homothetic).

\begin{remar}\rm
If $\Omega$ is not irreducible, then its painted Dynkin diagram is given by the disjoint union of (connected) painted Dynkin diagrams of simple groups, and the coordinates are Bochner if and only if are the coordinates on each factor.
\end{remar}

\begin{remar}\rm\label{remarteor}
Theorem \ref{thm main} represents only a partial answer to Question 1  for two reasons. First, we are assuming that the flag manifolds are of classical type. Secondly, we are showing that $\Boc_0=V_0$ only for those flag manifolds (of classical type) satisfying conditions (i), (ii) and (iii).
We believe that  the equality  $\Boc_0=V_0$  still holds  also for the all  flag manifolds of classical type, but we were not able to  attack this problem due to the difficulty of the  computations involved.
\end{remar}

The paper is organized as follows. In Section 2, we recall the basic facts on flag manifolds and we provide the details of the construction given in \cite{A-P}. Section 3 is dedicated to the proof of Theorem \ref{thm main}.

\section{Complex coordinates on the flag manifold $\Omega = G/K$}

Let us begin by recalling some basic facts about the theory of flag manifolds and painted Dynkin diagrams.

Given a compact semisimple group $G$ and $Z \in \mathfrak{g}$, the orbit $F =\mathrm{Ad}_G Z$ of $Z$ for the adjoint action of $G$ on $\mathfrak{g}$ is diffeomorphic to the quotient manifold $G/K$, being $K$ the stabilizer of $Z$ with respect to the adjoint action, and is called a {\it flag manifold}. 

Fixed a Cartan subalgebra $\mathfrak{h}^{\bC}$ of the complexification $\mathfrak{g}^{\bC}$ of $\mathfrak{g}$, a {\it root} $\alpha \in R$ is a functional $\mathfrak{h}^{\bC} \rightarrow \bC$ such that $[H, E_{\alpha}] = \alpha(H) E_{\alpha}$ for each $H \in \mathfrak{h}^{\bC}$ and for some $E_{\alpha} \in \mathfrak{g}^{\bC}$ (called the root vector of $\alpha$)

Up to replacing $Z$ with another point of the orbit $\mathrm{Ad}_G Z$, one can assume that $Z$ belongs to $\mathfrak{h}^{\bC}$ and denote $R_K = \{ \alpha \in R \ | \ \alpha(Z) = 0 \}$. Then, one has the decomposition

\begin{equation}\label{decompositionK}
\mathfrak{k}^{\bC} = \mathfrak{h}^{\bC} + \sum_{\alpha \in R_K} \bC E_{\alpha} .
\end{equation}

The set $R_K$ can be in fact also described as the root system of the semisimple part of $\mathfrak{k}^{\bC}$ and its elements (resp. the elements of the complementary subset $R_M$) are usually called {\it white roots} (resp. {\it black roots}). The reason is that one can represent the flag manifold $G/K$ on the Dynkin diagram of $G$, equipped with a given basis\footnote{A {\it basis} $\Pi$ of the root system is a subset $\Pi \subseteq R$ such that every root $\alpha \in R$ can be written as a linear combination of the elements of $\Pi$ with the coefficients either all non-negative or all non-positive. In the first (resp. second) case, $\alpha$ is said to be positive (resp. negative). The set of positive roots with respect to a fixed basis will be denoted by $R^+$.} $\Pi$, by painting black the vertices corresponding to roots belonging to $R_M$. One gets a decomposition $\Pi = \Pi_K \cup \Pi_M$ of the basis $\Pi$ and the resulting diagram is called {\it painted Dynkin diagram}. Looking at a painted diagram, one can easily recover the root decomposition and the flag manifold: indeed, a root $\alpha$ belongs to $R_M$ if and only if $\alpha = \sum_{\beta \in \Pi} c_{\beta} \beta$ with $c_{\beta} \neq 0$ for some $\beta \in \Pi_M$; moreover, deleting the black nodes from the diagram one gets the Dynkin diagram of the semisimple part of $K$.

The number of black nodes in the painted diagram equals the second Betti number $b_2(G/K)$. This can be seen, for example, by deriving both hand-sides of the following formula, proved in \cite{Ak} (see Corollary 3.2), which gives an explicit expression for the Poincar\'e polynomial of $G/K$:

\begin{equation}\label{formula betti}
P(G/K, t^{1/2}) = \prod_{\alpha \in R_{\mathfrak{m}}^+} \frac{1 - t^{h(\alpha)+1}}{1 - t^{h(\alpha)}}.
\end{equation}
where, for every root $\alpha$, the {\it height} $h(\alpha)$ is defined as the sum $h(\alpha )= \sum_{i=1}^m k_i$, being $\alpha = \sum_{i=1}^m k_i \alpha_i \in R^+$ the decomposition of $\alpha$ in terms of the basis $\Pi = \{\alpha_1, \dots, \alpha_m \}$.
\begin{Ex}\label{example classicals}(\cite{helgason}, Section III.8)
\rm
\begin{enumerate}
\item [$G= SU(n)$ : ] $\mathfrak{g}^{\bC} = sl(n, \bC)$ is the set of matrices with null trace, and a Cartan subalgebra $\mathfrak{h}^{\bC}$ is given by the diagonal matrices in $sl(n, \bC)$; for any $H = \diag(h_1, \dots, h_n)$ let $e_i(H) = h_i$: then the root system is $R = \{ e_i - e_j \ | \ i \neq j \}$ and $E_{\alpha}$, $\alpha = e_i - e_j$, is the matrix $E_{ij}$ having 1 in the ij place and 0 anywhere else. 
Given the basis $\Pi_{can} = \{ \alpha_1 = e_1 - e_2, \dots, \alpha_{n-1} = e_{n-1} - e_n \}$, the Dynkin diagram  is
\begin{align*}
&&& \underset{\substack{\alpha_1}}{\circ} - \underset{\substack{\alpha_2}}{\circ} - \dotsb - \underset{\substack{\alpha_{n-2}}}{\circ} - \underset{\substack{\alpha_{n-1}}}{\circ} &&
 \\
\end{align*}
\item [$G= Sp(n)$ : ]$\mathfrak{g}^{\bC} = sp(n, \bC)$ is the set of $2n \times 2n$ block matrices of the kind $\left( \begin{array}{cc}
Z_1 & Z_2 \\
Z_3 & -{}^T Z_1 \\
\end{array} \right)$, where $Z_2, Z_3$ are symmetric. A Cartan subalgebra $\mathfrak{h}^{\bC}$ is given by diagonal matrices 
$$H = \diag(h_1, \dots, h_n, -h_1, \dots, -h_n)$$
in $sp(n, \bC)$, and if for any such $H$ we define $e_i(H) = h_i$, $i=1, \dots, n$, then the root system is $R = \{ \pm e_i \pm e_j \}$ (the case $i=j$ is allowed when the signs are equal). The root vector $E_{\alpha}$ is given by
$\left( \begin{array}{cc}
E_{ij} & 0 \\
0 & -E_{ji} \\
\end{array} \right)$ if $\alpha = e_i - e_j$, $\left( \begin{array}{cc}
0 & E_{ij} + E_{ji} \\
0 & 0 \\
\end{array} \right)$ if $\alpha = e_i + e_j$ and $\left( \begin{array}{cc}
0 & 0 \\
E_{ij} + E_{ji} & 0 \\
\end{array} \right)$ if $\alpha = -e_i - e_j$.
Given the basis $\Pi_{can} = \{ \alpha_1 = e_1 - e_2, \dots, \alpha_{n-1} = e_{n-1} - e_n, \alpha_n = 2 e_n \}$, the Dynkin diagram is 
\begin{align*}
&&& \underset{\substack{\alpha_1}}{\circ} - \underset{\substack{\alpha_2}}{\circ} - \dotsb - \underset{\substack{\alpha_{n-1}}}{\circ} \Leftarrow \underset{\substack{\alpha_n}}{\circ} && \\
\end{align*}
\item [$G= SO(2k)$ : ] here and throughout the paper we identify the complexification $SO(2k, \bC)$ with the subgroup of $GL(2k, \bC)$ leaving invariant the quadratic form $z_1 z_{k+1} + \cdots + z_k z_{2k}$. Then $\mathfrak{g}^{\bC} = so(2k, \bC)$ is the set of $2k \times 2k$ block matrices of the kind $\left( \begin{array}{cc}
Z_1 & Z_2 \\
Z_3 & -{}^T Z_1 \\
\end{array} \right)$, where $Z_2, Z_3$ are skew-symmetric. A Cartan subalgebra $\mathfrak{h}^{\bC}$ is given by diagonal matrices $H = \diag(h_1, \dots, h_k, -h_1, \dots, -h_k)$ in $so(2k, \bC)$, and if for any such $H$ we define $e_i(H) = h_i$, $i=1, \dots, k$, then the root system is $R = \{ \pm e_i \pm e_j \ (i \neq j) \}$. The root vector $E_{\alpha}$ is given by
$\left( \begin{array}{cc}
E_{ij} & 0 \\
0 & -E_{ji} \\
\end{array} \right)$ if $\alpha = e_i - e_j$, $\left( \begin{array}{cc}
0 & E_{ij} - E_{ji} \\
0 & 0 \\
\end{array} \right)$ if $\alpha = e_i + e_j$ ($i<j$) and $\left( \begin{array}{cc}
0 & 0 \\
E_{ij} - E_{ji} & 0 \\
\end{array} \right)$ if $\alpha = -e_i - e_j$ ($i<j$).
Given the basis $\Pi_{can} = \{ \alpha_1 = e_1 - e_2, \dots, \alpha_{k-1} = e_{k-1} - e_k, \alpha_k = e_{k-1} + e_k \}$, the Dynkin diagram is
\begin{align*}
&&& \underset{\substack{\alpha_1}}{\circ} - \underset{\substack{\alpha_2}}{\circ} - \dotsb - \underset{\substack{\alpha_{k-2}}}{\overset{\overset{\textstyle\circ_{\alpha_k}}{\textstyle\vert}}{\circ}} \,-\, \underset{\substack{\alpha_{k-1}}}{\circ} && \\
\end{align*}
\item [$G= SO(2k+1)$ : ] here and throughout the paper we identify the complexification $SO(2k+1, \bC)$ with the subgroup of $GL(2k+1, \bC)$ leaving invariant the quadratic form $2(z_1 z_{k+1} + \cdots + z_k z_{2k}) + z_{2k+1}$. Then $\mathfrak{g}^{\bC} = so(2k+1, \bC)$ is the set of $(2k+1) \times (2k+1)$ block matrices of the kind $\left( \begin{array}{ccc}
Z_1 & Z_2 & u \\
Z_3 & -{}^T Z_1 & v \\
-{}^T v & -{}^T u & 0
\end{array} \right)$, where $Z_2, Z_3$ are skew-symmetric and $u, v \in \bC^k$. A Cartan subalgebra $\mathfrak{h}^{\bC}$ is given by diagonal matrices $H = \diag(h_1, \dots, h_k, -h_1, \dots, -h_k, 0)$ in $so(2k+1, \bC)$, and if for any such $H$ we define $e_i(H) = h_i$, $i=1, \dots, k$, then the root system is $R = \{ \pm e_i \pm e_j \ (i \neq j), \pm e_i \}$. The root vector $E_{\alpha}$ is given by
$\left( \begin{array}{ccc}
E_{ij} & 0 & 0 \\
0 & -E_{ji} & 0\\
0 & 0 & 0 
\end{array} \right)$ if $\alpha = e_i - e_j$, $\left( \begin{array}{ccc}
0 & E_{ij} - E_{ji} & 0 \\
0 & 0 & 0 \\
0 & 0 & 0
\end{array} \right)$ if $\alpha = e_i + e_j$ ($i<j$), $\left( \begin{array}{ccc}
0 & 0 & 0 \\
E_{ij} - E_{ji} & 0 & 0 \\
0 & 0 & 0
\end{array} \right)$ if $\alpha = -e_i - e_j$ ($i<j$), $\left( \begin{array}{ccc}
0 & 0 & E_i \\
0 & 0 & 0 \\
0 & -{}^T E_i & 0
\end{array} \right)$ if $\alpha = e_i $ and $\left( \begin{array}{ccc}
0 & 0 & 0 \\
0 & 0 & E_i \\
-{}^T E_i & 0  & 0
\end{array} \right)$ if $\alpha = - e_i $ (being $E_i$ the $i$-th vector of the canonical basis of $\bC^k$).

Given the basis $\Pi_{can} = \{ \alpha_1 = e_1 - e_2, \dots, \alpha_{k-1} = e_{k-1} - e_k,  \alpha_k = e_k \}$, the Dynkin diagram is
\begin{align*}
&&& \underset{\substack{\alpha_1}}{\circ} - \underset{\substack{\alpha_2}}{\circ} - \dotsb - \underset{\substack{\alpha_{k-1}}}{\circ} \Rightarrow \underset{\substack{\alpha_k}}{\circ} && \\
\end{align*}
\end{enumerate}
\end{Ex}

Now, it is known that invariant complex structures on a flag manifold $G/K$ are in one-to-one correspondence with {\it maximal closed nonsymmetric subsets} $Q$ of the set $R_M$ of the black roots of $G/K$:

\begin{Def}\label{defQ}
A subset $Q \subseteq R_M$ is said to be maximal closed nonsymmetric if it satisfies the following conditions:

\begin{enumerate}

\item[(i)] $Q \cup (-Q) = R_M$;

\item[(ii)] $Q \cap (-Q) = \emptyset$;

\item[(iii)] for any $\alpha, \beta \in Q$ such that $\alpha + \beta \in R$ one has $\alpha + \beta \in Q$.

\end{enumerate}

\end{Def}

More precisely, the manifold $G/K$ endowed with the complex structure $J_Q$ corresponding to $Q$ is biholomorphic to the complex homogeneous manifold $G^{\bC}/K^{\bC} G^{Q}$, where $G^{Q} = \exp(\mathfrak{g}^{Q})$ and $\mathfrak{g}^{Q} = \sum_{\alpha \in Q} \bC E_{\alpha}$. 

Since the product $G_{reg}^{\bC} = G^{-Q} K^{\bC} G^{Q}$ (where $G^{-Q} = \exp(\mathfrak{g}^{-Q})$ and $\mathfrak{g}^{-Q} = \sum_{\alpha \in -Q} \bC E_{\alpha}$) defines an open dense subset in $G^{\bC}$, its image in $G^{\bC}/K^{\bC} G^{Q}$ via the natural projection $G^{\bC} \rightarrow G^{\bC}/K^{\bC} G^{Q}$ defines an open dense subset in $G/K$, denoted $F_{reg} = G_{reg}^{\bC}/K^{\bC} G^Q$. Clearly, $F_{reg} \simeq G^{-Q}$.
 Then, by

\begin{equation}\label{complexcoordinates}
z = (z_{\alpha})_{\alpha \in -Q} \in \bC^N \mapsto \exp(Z(z)) \in G^{-Q} \simeq F_{reg} \subseteq F
\end{equation}
where
\begin{equation}\label{Z(z)}
Z(z) = \sum_{\alpha \in -Q} z_{\alpha} E_{\alpha}
\end{equation}
(where $N$ is the cardinality of $Q$) one defines a system of complex coordinates on $F_{reg}$. 

\begin{remar}\label{remarkalternativeQ}
\rm Given a flag manifold $F = G/K$ endowed with an invariant complex structure, then, up to $G$-diffeomorphism, one can always assume that $F$ is represented by a painted Dynkin diagram endowed with the canonical equipment $\Pi_{can}$ (Examples \ref{example classicals}) and that the complex structure is associated to $Q = R_M^+ := R_M \cap R^+$ (the positive roots with respect to $\Pi_{can}$). For more details, see \cite{LMZ}.
\end{remar}
Now we show how to find an explicit K\"ahler potential of any invariant K\"ahler form $\omega$  in the coordinates (\ref{complexcoordinates}) in the case when $G$ is one of the classical groups $SU(n), Sp(n), SO(n)$ (where $SO(n)$ is realized as a group of matrices as in Examples \ref{example classicals}).

Let us first recall that, given the decomposition $\Pi = \Pi_K \cup \Pi_M$ of the basis $\Pi$ into white and black roots, where $\Pi_K = \{ \beta_1, \dots, \beta_k\}$ and $\Pi_M = \{ \a_1, \dots, \a_m\}$, then the {\it fundamental weight} $\bar \a_i$
associated with $\a_i$, $i = 1, \dots, m$ is the element of ${\mathfrak{h}}^*$ defined by
\begin{equation}\label{deffundwe}
\frac{2 \langle \bar \alpha_i, \alpha_j \rangle}{\| \alpha_j \|^2} = \delta_{ij}, \ \ \langle \bar \alpha_i, \beta_j \rangle = 0
\end{equation}
where $\langle, \rangle$ denotes the scalar product induced on the real space ${\mathfrak{h}}^*$ spanned by the roots by the Killing form of $G$.
If we denote by $\mathfrak{t} = Z(\mathfrak{k}) \cap \mathfrak{h}$ the intersection between the center $Z(\mathfrak{k})$ of $\mathfrak{k}$ and $\mathfrak{h}$ , then the fundamental weights form a basis of the real space $\mathfrak{t}^*$.

\begin{Prop}\label{transgression}(\cite{A-P}, Proposition 2.2) There exists a natural isomorphism $\xi \mapsto \omega_{\xi}$  between $\mathfrak{t}^*$ and the space of invariant 2-forms on $F$; moreover, if $J_Q$ is the complex structure associated to $Q = R_M^+$ with the given equipment, then $\omega_{\xi}$ is K\"ahler with respect to $J_Q$ if and only if all the coordinates of $\xi$ with respect to the fundamental weights $\bar \alpha_1, \dots, \bar \alpha_m$ are positive.
\end{Prop}
In order to define a potential for the \K\ metric $\omega_{\xi}$, we need to give the following

\begin{Def}\label{admissible}(\cite{A-P}, Definition 8.1)
Let $F = G/K$, $G \subseteq GL(n, \bC)$, be a flag manifold. A principal minor $\Delta_k$, $k=1, \dots, n-1$, (i.e. the function $GL(n, \bC) \rightarrow \bC$ associating to $A \in GL(n, \bC)$ the determinant of the submatrix $A_k$ of $A$ given by the first $k$ rows and columns of $A$) is said to be $F$-{\it admissible} if for every $A \in K^{\bC}$ and every $v = (v_1, \dots, v_n) \in \bC^n$, $v_{k+1} = \cdots = v_n=0$ implies $(vA)_{k+1} = \cdots = (vA)_n = 0$.

\end{Def}

\begin{Ex}\label{exsflags}\rm
For the flag manifolds of the classical groups (see, for example, \cite{arv}) 
\begin{enumerate}

\item[] $G/K = SU(n)/S(U(n_1) \times \cdots \times U(n_s))$ ($n = n_1 + \cdots + n_s$, $s \geq 1$): 

\item[] $G/K = Sp(n)/U(n_1) \times \cdots \times U(n_s) \times Sp(l)$

\item[] $G/K = SO(2n+1)/U(n_1) \times \cdots \times U(n_s) \times SO(2l+1)$

\item[] $G/K = SO(2n)/U(n_1) \times \cdots \times U(n_s) \times SO(2l)$

($n = n_1 + \cdots + n_s + l$, $s, l \geq 0, l \neq 1$)

\end{enumerate}
it is easy to see that a minor $\Delta_l$ is admissible if and only if $l = n_1 + \cdots + n_j$, for some $j=1, \dots, s-1$ (resp. $j=1, \dots, s$) in the case $G=SU(n)$ (resp. in all the other cases).
\end{Ex}

We have the following
\begin{Th}\label{potential}(\cite{A-P}, Proposition 8.2)
Let $F = G/K$ be a flag manifold of classical type represented by a painted Dynkin diagram endowed with the canonical equipment $\{ \alpha_1, \dots, \alpha_n \}$ and let $\{ \alpha_{j_1}, \dots, \alpha_{j_s} \}$, $j_1 < \cdots < j_s$, be the set of black nodes. Let $F$ be endowed with the $G$-invariant complex structure determined by $Q = R_M^+$. 
Then, in the holomorphic coordinates $z = (z_{\alpha})_{\alpha \in -Q}$ defined in (\ref{complexcoordinates}), a K\"ahler potential for the K\"ahler form $\omega_{\xi}$, $\xi = \sum_{k=1}^s c_k \bar \alpha_{j_k}$, with $c_1, \dots, c_s > 0$  is 

\begin{equation}\label{potential}
D_0(z) = \sum_{k=1}^s c_k \ln \Delta_{l_k}({}^T \overline{\exp(Z(z))} \exp(Z(z)) )
\end{equation}
where  $\Delta_{l_1}, \dots, \Delta_{l_s}$ are admissible minors.
\end{Th}
Notice that, by Remark \ref{remarkalternativeQ}, the assumptions that $\{ \alpha_1, \dots, \alpha_n \}$ is the canonical equipment and that  $Q = R_M^+$ are not restrictive.

\section{Proof of Theorem \ref{thm main}}
We have already pointed out in the introduction the K\"ahler potential (\ref{potential}) coincides with Calabi's diastasis function centered at the origin.
Therefore, by definition,  in order to see when (\ref{complexcoordinates}) are Bochner we must check if in the expansion of the diastasis (\ref{potential}) there are terms of the kind $z_{\alpha_1} \bar z_{\alpha_2} \cdots \bar z_{\alpha_k}$ (or their conjugates): from now on, we will call these terms {\it forbidden monomials}.
As we will see below, it will be sufficient to check the existence of forbidden monomials of degree 3, i.e. forbidden trinomials.
Since (\ref{potential}) is the diastasis centered in the origin, the monomials of its expansion have holomorphic and anti-holomorphic part both non trivial: this fact together with $\ln(1+x)= x + O(x^2)$, implies that the trinomials of the expansion of  (\ref{potential}) are exactly the ones of 

\begin{equation}\label{diastasisSENZAlog}
\sum_{k=1}^{s} c_k \,\Delta_{l_k}\left({}^T \overline{\exp\left(Z\left(z\right)\right)} \exp\left(Z\left(z\right)\right) \right).
\end{equation}

In order to determine these trinomials, we prove the following, technical

\begin{Lem}\label{lem trinomial}
Assume that $Z$ is a $m\times m$ matrix, and $r \leq m$; then the trinomials of type $\ov Z_{i_1j_1}\, \ov Z_{i_2j_2} \,  Z_{i_3j_3}$ of
 $$\Delta_{r}\left({}^T \overline{\exp\left(Z\right)} \exp\left(Z\right) \right)$$ 
are of the following types:

\begin{enumerate}
\item[(I)] $+\frac{1}{2} Z_{si} \ov Z_{st} \ov Z_{ti}$, $i \leq r$ and $s, t= 1, \dots, m$
\item[(II)] $-\frac{1}{2} Z_{ij} \ov Z_{is} \ov Z_{sj}$, $i, j \leq r$, $i \neq j$ and $s= 1, \dots, m$
\item[(III)] $-Z_{si} \ov Z_{sj} \ov Z_{ji}$, $i, j \leq r$, $i \neq j$ and $s= 1, \dots, m$
\item[(IV)] $+ Z_{ab} \ov Z_{ac} \ov Z_{cb}$, $\{ a, b, c \} = \{ i, j, k \}$, for $i,j,k \leq r$, $i<j<k$
\end{enumerate}
\end{Lem}

\proof
Since  $Z^k$ is a matrix whose entries are homogeneous polynomial of degree $k$,
it is enough to study the terms of $\Delta_{r}A$, where
 
$$A := {}^T \overline{\left(I+ Z + \frac{1}{2} Z^2\right)}\, \left(I+ Z + \frac{1}{2} Z^2\right)=$$
$$= I + Z + {}^T \ov Z + {}^T \ov Z Z + \frac{1}{2} {}^T \ov Z^2 + \frac{1}{2} {}^T \ov Z^2 Z + \cdots$$ 

By the very definition of the determinant  and by the assumption on $Z$ we see that the terms of third degree of $\Delta_r(Z,\ov Z)$ come from addendums of the form
$$
A_{11}\cdots \widehat A_{ii} \cdots \widehat A_{jj} \cdots \widehat A_{kk} \cdots A_{rr}\cdot A_{i\sigma(i)}\cdot A_{j\sigma(j)}\cdot A_{k\sigma(k)}.
$$ 
that is we have to consider the three cases 
\begin{enumerate}
\item[(a)] +$A_{ii}$ for every $i \leq r$; 
\item[(b)] -$A_{ij}A_{ji}$ for every $i, j \leq r$, $i < j$; 
\item[(c)] +$A_{ij}A_{jk}A_{ki}$ and $A_{ik}A_{ji}A_{kj}$ for every $i,j,k \leq r$, $i<j<k$.
\end{enumerate}

In case (a), by the definition of $A$ we have that the trinomials of the statement are given by
\begin{equation}\label{Aii}
+\frac{1}{2} ({}^T \ov Z^2 Z)_{ii} = +\frac{1}{2} \ov Z_{st} \ov Z_{ti} Z_{si}, \ \ s,t =1, \dots, d 
\end{equation}
which corresponds to case (I) of the statement.
In case (b), the trinomials of the statement are given by
\begin{equation}\label{Aij2}
-Z_{ij} \cdot \frac{1}{2} ({}^T \ov Z^2)_{ji} = -\frac{1}{2} Z_{ij} \ov Z_{is} \ov Z_{sj}, \ \ -\frac{1}{2} ({}^T \ov Z^2)_{ij} \cdot Z_{ji}  = -\frac{1}{2} Z_{ji} \ov Z_{js} \ov Z_{si}
\end{equation}
\begin{equation}\label{Aij3}
-{}^T \ov Z_{ij} \cdot ({}^T \ov Z Z)_{ji} = - \ov Z_{ji} \ov Z_{sj} Z_{si}, \ \ -({}^T \ov Z Z)_{ji} \cdot {}^T \ov Z_{ij}  = - \ov Z_{ij} \ov Z_{si} Z_{sj}
\end{equation}
for $s =1, \dots, d$, and correspond to cases (II) and (III) of the statement.
In case (c), for $A_{ij}A_{jk}A_{ki}$ we get the trinomials 
\begin{equation}\label{Aijk1}
Z_{ij} \ov Z_{kj} \ov Z_{ik}, \ \ \ov Z_{ji} Z_{jk} \ov Z_{ik}, \ \ \ov Z_{ji} \ov Z_{kj} Z_{ki}
\end{equation}
while for $A_{ik}A_{ji}A_{kj}$ we get the trinomials  
\begin{equation}\label{Aijk2}
Z_{ik} \ov Z_{ij} \ov Z_{jk}, \ \ \ov Z_{ki} Z_{ji} \ov Z_{jk}, \ \ \ov Z_{ki} \ov Z_{ij} Z_{kj}
\end{equation}
which give the case (IV) of the statement. The proof is done.
\endproof

Now we are ready to prove Theorem \ref{thm main}. 

Recall that, as seen in the previous section, the second Betti number $b_2(\Omega)$ of $\Omega= G/K$ equals the number of black nodes in the painted diagram of $G/K$. So, we are going to  distinguish the cases where the painted diagram of $\Omega$ has ony one black node and the case where there are at least two black nodes.

For the sake of simplicity, from now on, by {\it trinomial of type (I)-(IV)} we will refer to the types determined in the statement of Lemma \ref{lem trinomial}; moreover, the matriz $Z(z)$ defined in (\ref{Z(z)}) will be denoted by $Z$. Clearly, $Z_{ij}$ is not identically zero if and only if $Z_{ij} = \pm z_{\alpha}$ for some $\alpha \in -Q$.

\subsection{The case of second Betti number $b_2(\Omega) =1$} 
In this case, the painted Dynkin diagram of $\Omega$ has only one black root, which can be either $e_r-e_{r+1}$ (for any group) or $e_{d-1} + e_d, 2 e_d, e_d$ (respectively for $G = SO(2d), Sp(d), SO(2d+1)$). In all the cases, we are going to prove that the coordinates are Bochner, that is (i) of Theorem \ref{thm main}.

\subsubsection{$G = SU(d)$, $\Pi_M = \{e_r-e_{r+1} \}$} \hfill

In this case, we have
$$
Z=\begin{pmatrix} 0 & 0 \\ * & 0\end{pmatrix},
$$
and it is immediately seen that $Z^2 = 0$, so that $\exp\left(Z \right)=I+Z$ and, since for $i,j\leq r$ we have $Z_{ij}=0$,
the entry $({}^T \overline{\exp\left(Z \right)} \exp\left(Z \right))_{ij}$ is given by
$$\sum_{k=1}^{d} \ov{ \left(\delta_{ki} + Z_{ki}\right)} \left(\delta_{kj} + Z_{kj}\right) = \delta_{ij} +\sum_{k=1}^{d} \ov Z_{ki}\, Z_{kj}.
$$
Therefore the (only) admissible minor $\Delta_{r}$ is given by
$$
\Delta_{r}=1 + \left\| z \right\|^2 +  \psi_{2,2},
$$
where $z = (z_{\alpha})_{\alpha \in -Q}$ and $\psi_{2,2}$ is a power series whose terms have degree $\geq 2$ in both the variables $z_\alpha$ and $\bar z_\alpha$, which implies that \eqref{complexcoordinates} are Bochner.  The proof is complete. 

\subsubsection{$G = Sp(d), SO(2d)$, $\Pi_M = \{e_r-e_{r+1} \}$} \hfill
The only admissible minor is $\Delta_r$. On the one hand, since $Z_{kl} = 0$ if $k, l \leq r$, we see that there are no trinomials of type (II),(III),(IV); on the other hand, there are no trinomials $+\frac{1}{2} Z_{si} \ov Z_{st} \ov Z_{ti}$ of type (I) with all the indices less or equal to $d$, since this would imply that $e_i - e_t, e_t - e_s \in R_M^+$, that is, being $e_r-e_{r+1}$ the only black node, $t \geq r+1$ and $t \leq r$, a contradiction. 

Moreover, $Z_{st} = 0$ if $s \leq d, t \geq d+1$, so we are left with the following forbidden trinomials

$$\frac{1}{2} \sum_{i=1}^r \sum_{s,t=d+1}^{2d} Z_{si} \ov Z_{st} \ov Z_{ti} + \frac{1}{2} \sum_{i=1}^r \sum_{s=d+1}^{2d} \sum_{t=1}^{d} Z_{si} \ov Z_{st} \ov Z_{ti}$$

Now we show that this expression is in fact zero by showing that the two sums simpify.
Indeed, in both cases $G=SO(2d)$ and $G=Sp(d)$, by the symmetries of the Lie algebra, by denoting $s = \tilde s + d$ and $t = \tilde t + d$, the first sum writes 

$$\frac{1}{2} \sum_{i=1}^r \sum_{s,t=d+1}^{2d} Z_{si} \ov Z_{st} \ov Z_{ti} = - \frac{1}{2} \sum_{i=1}^r \sum_{\tilde t=1}^{d} \sum_{\tilde s=1}^{r}  Z_{i+d, \tilde s} \ov Z_{\tilde t \tilde s} \ov Z_{i+d, \tilde t} = $$

\begin{equation}\label{sum1radice}
- \frac{1}{2} \sum_{\tilde i=d+1}^{d+r} \sum_{\tilde t=1}^{d} \sum_{\tilde s=1}^{r}  Z_{\tilde i \tilde s} \ov Z_{\tilde t \tilde s} \ov Z_{\tilde i \tilde t}
\end{equation}
(where $\tilde s \leq r$ since $Z_{\tilde t \tilde s} =0$ for $\tilde t \leq d$ and $\tilde s > r$).
As for the second sum, we notice that if $s > r + d$ then it must be $t \leq r$ (otherwise $Z_{st} = 0$), but then we would have $Z_{ti} = 0$. So it rewrites
$$\frac{1}{2} \sum_{i=1}^r \sum_{s=d+1}^{d+r} \sum_{t=1}^{d} Z_{si} \ov Z_{st} \ov Z_{ti}$$
and it simplifies with (\ref{sum1radice}), as claimed.
Since it is easily seen that $Z^3 = 0$ (use the fact that $Z^2 \neq 0$ only if $d < l \leq d+r$ and $1 \leq t \leq r$), there cannot be other forbidden monomials, so the coordinates are Bochner.

\subsubsection{$G = SO(2d+1)$, $\Pi_M = \{e_r-e_{r+1} \}$} \hfill
We can use the same arguments as in the above case $G = SO(2d)$ to see that there are no forbidden trinomials with indices less or equal to $2d$: so, we are left with
$$\frac{1}{2} \sum_{i=1}^r \sum_{t=1}^{2d} Z_{2d+1, i} \ov Z_{2d+1, t} \ov Z_{ti} + \frac{1}{2} \sum_{i=1}^r \sum_{s=1}^{2d}  Z_{si} \ov Z_{s, 2d+1} \ov Z_{2d+1, i}.$$
But the first sum vanishes since, in order to have $Z_{2d+1, t} \neq 0$ it must be $t \leq r$, which gives $Z_{ti}=0$; the second sum vanishes because in order to have $Z_{s, 2d+1}\neq 0$ it must be $s = d+1, \dots, d+r$, and then by the symmetries of $SO(2d+1)$ it rewrites

$$\sum_{i=1}^r \sum_{\tilde s=1}^{r}  Z_{\tilde s+d,i} \ov Z_{\tilde s+d, 2d+1} \ov Z_{2d+1, i} = - \sum_{i=1}^r \sum_{\tilde s=1}^{r}  Z_{i+d,\tilde s} \ov Z_{2d+1, \tilde s} \ov Z_{i+d, 2d+1}$$
and the second equality means that it is zero, as claimed.
As above, one easily proves that $Z^3=0$, which show  that also in this case the coordinates are Bochner.

\subsubsection{$Sp(d)$, $\Pi_M = \{2 e_{d}\}$ and $SO(2d)$, $\Pi_M = \{e_{d-1} + e_{d} \}$} \hfill
Since in both cases we have
$$
Z=\begin{pmatrix} 0 & 0 \\ * & 0\end{pmatrix},
$$
we can use the same arguments as in the above case $G = SU(d)$, $\Pi_M = \{e_r-e_{r+1} \}$ which prove that the coordinates are Bochner.

\subsubsection{$SO(2d+1)$, $\Pi_M = \{ e_{d} \}$} 
In this case
$$
Z=\begin{pmatrix} 0 & 0 & 0 \\ * & 0 & u \\ -{}^T u & 0 & 0 \end{pmatrix},
$$
from which  one easily sees by a straight computation that $Z^2 \neq 0$ (and $Z^3=0$). 
In order to see that there are no forbidden trinomials, notice that since $Z_{kl} = 0$ if $k,l \leq d$, there are no trinomials of type (II),(III),(IV), and for the same reason in the trinomials $+\frac{1}{2} Z_{si} \ov Z_{st} \ov Z_{ti}$ of type (I) it must be $t, s = d+1, \dots, 2d+1$: this in turn implies that in fact $t = 2d+1$ in order to have $Z_{st} \neq 0$. So we are left with

$$+\frac{1}{2} \sum_{i=1}^{d} \sum_{s=d+1}^{2d} Z_{si} \ov Z_{s, 2d+1} \ov Z_{2d+1, i} = \frac{1}{2} \sum_{i=1}^{d} \sum_{\tilde s=1}^{d} Z_{\tilde s + d, i} \ov Z_{\tilde s + d, 2d+1} \ov Z_{2d+1, i}$$
which, by the symmetries of $SO(2d+1)$ equals 
$$-\frac{1}{2} \sum_{i=1}^{d} \sum_{\tilde s=1}^{d} Z_{i + d, \tilde s} \ov Z_{2d+1, \tilde s} \ov Z_{i + d, 2d+1},$$ 
and hence it is forced to be  zero. So  the coordinates are Bochner.

\subsection{The case of second Betti number $b_2(\Omega) \geq 2$} 
The assumption means that the painted Dynkin diagram of $\Omega$ has at least two black nodes. 
We first consider the two Bochner cases stated in (ii) and (iii) of Theorem \ref{thm main}: 

\subsubsection{$G = SU(d)$, $\Pi_M$ contains exactly two nodes (necessarily of type $e_j-e_{j+1}$)} 
Let us assume that $e_k-e_{k+1}, e_r-e_{r+1}$ are the black nodes, so that $\Delta_k$ and $\Delta_r$ are the admissible minors. Let $c_k, c_r$ be the coefficients in front of these trinomials in the expansion of the potential.
First, we notice that $Z^3=0$, since for any $i,j$ we have 
$$Z^3_{ij} = \sum_{k,l=1}^n Z_{ik} Z_{kl} Z_{lj},$$ 
so that in order to have a non zero term in this sum we should take indices $j<l<k<i$ such that $e_k - e_i, e_l-e_k, e_j-e_l \in R_M^+$, and it is easy to see that this is not possible with only two black nodes.
Then, we have to check just forbidden trinomials.
Now, in $\Delta_k$ we have only trinomials $+\frac{1}{2} Z_{si} \ov Z_{st} \ov Z_{ti}$, $i \leq k$ and $s, t= 1, \dots, d$ of type (I) since if $i, j \leq k$ then $Z_{ij} = 0$ because there are no black nodes before $e_k-e_{k+1}$; moreover, by the structure of the Lie algebra and the fact that we have only those black nodes, we must have that $t$ and $s$ must be both greater than $k$, and more precisely $s > r$ and $k < t \leq r$.
On the other hand, in $\Delta_r$ we have trinomials of all the four types of the statement of Lemma \ref{lem trinomial}: more precisely, monomials of type (II) simplify with monomials of type (I) with $s \leq r$, so that from these two types we are left with monomials of the type 

\begin{equation}\label{1e2}
+\frac{1}{2} Z_{si} \ov Z_{st} \ov Z_{ti}, \ \  i \leq r, s > r, t \leq r
\end{equation} 
(the last condition follows from $s > r$ in order to have non zero $Z_{st}$). 

Similarly, in $\Delta_r$ monomials of type (IV) simplify with monomials of type (III) with $s \leq r$, so from these types of monomials we are left with 
\begin{equation}\label{3e4}
-Z_{si} \ov Z_{sj} \ov Z_{ji}, \ \ i, j \leq r, s > r
\end{equation}
Summing (\ref{1e2}) and (\ref{3e4}) (rename $t$ as $j$ in (\ref{1e2}) ) we finally see that in $\Delta_r$ we have just
\begin{equation}\label{finallyDeltar}
-\frac{1}{2} Z_{si} \ov Z_{sj} \ov Z_{ji}, \ \ i, j \leq r, s > r
\end{equation}
In fact, since there are only two black roots we must have $i \leq k$ and $k < j \leq r$.
Then, comparing with the trinomials in $\Delta_k$ we see that they appear with opposite sign, so that in the expansion they simplify if and only if $c_k = c_r$.

\subsubsection{$SO(2d)$, $\Pi_M = \{ e_{1} - e_{2}, e_{d-1} + e_{d}\}$}
First of all, we notice that $Z^3=0$. Indeed, it is easy to see from the very structure of the root set in this case that $Z^3_{ij} = 0$ for $(i,j) \neq (d+1, 1)$; as for $Z^3_{d+1,1}$, we have
$$Z^3_{d+1,1} = \sum_{k=1}^{2d} Z_{d+1,k} Z^2_{k1} =$$
(since $Z^2_{k1} = \sum Z_{kl} Z_{l1} = Z_{k1} Z_{11} = 0$ for $k \leq d$)
$$= \sum_{k=d+1}^{2d} Z_{d+1,k} Z^2_{k1} = \sum_{j=1}^{d} Z_{d+1,d+j} Z^2_{d+j,1} = \sum_{j=1}^{d}  \sum_{i=1}^{2d} Z_{d+1,d+j} Z_{d+j,i}  Z_{i1}  = $$
(by the simmetries of $Z$ and the fact that $Z_{11} = 0$ and $Z_{d+j,i} = 0$ for $i > d$ and $j > 1$)
$$-  \sum_{j=1}^{d}  \sum_{i=1}^{2d} Z_{j1} Z_{d+j,i}  Z_{i1}  =-  \sum_{j=2}^{d}  \sum_{i=2}^{d} Z_{j1} Z_{d+j,i}  Z_{i1}.$$

This last expression is zero because, by the symmetries of $Z$ (i.e. $Z_{d+j,i}  = - Z_{d+i,j}$ ) and a change of indices we get
$$\sum_{j=2}^{d}  \sum_{i=2}^{d} Z_{j1} Z_{d+j,i}  Z_{i1} = - \sum_{j=2}^{d}  \sum_{i=2}^{d} Z_{j1} Z_{d+i,j}  Z_{i1} = - \sum_{i=2}^{d}  \sum_{j=2}^{d} Z_{i1} Z_{d+j,i}  Z_{j1}.$$

This implies that in order for the coordinates to be Bochner, we must just check that there are no fobidden trinomials.
In this case, the admissible minors are $\Delta_1$ and $\Delta_d$.
On the one hand, in $\Delta_d$ the trinomials of type (IV) vanish because, by the structure of the root system and the symmetries of the Lie algebra in this case, $Z_{kl} = 0$ provided either $k \leq d$ and $l > 1$ or $k=l=1$: then $Z_{ab} \ov Z_{ac} \ov Z_{cb}$ must have $b = c = 1$ and then $Z_{cb} = Z_{11} = 0$; as for trinomials $-\frac{1}{2} Z_{ij} \ov Z_{is} \ov Z_{sj}$ of type (II), since $i \leq d$ we must have $j=s=1$, but then $Z_{sj} = Z_{11} = 0$.
It remain only the trinomials of type (I) and (III).

As for type (I), that is  $\frac{1}{2}Z_{si} \ov Z_{st} \ov Z_{ti}$, $i \leq d$, and $s,t=1, \dots, 2d$, we notice that if $s \leq d$ then $i=t=1$ and then $Z_{ti} = Z_{11} = 0$; then we must have $s= d+1, \dots, 2d$. If $t=1, \dots, d$ then it must be $i=1$, that is we have $\frac{1}{2}Z_{s1} \ov Z_{st} \ov Z_{t1}$ ($s= d+1, \dots, 2d, t=1, \dots, d$); if $t= d+1, \dots, 2d$ then (always by the structure of the Lie algebra) it must be $s = d+1$, that is $\frac{1}{2}Z_{d+1 i} \ov Z_{d+1 t} \ov Z_{ti}$ ($t= d+1, \dots, 2d$), that is, by setting $t = \tilde t + d$, by the symmetries of the Lie algebra, $\frac{1}{2}Z_{d+1 i} \ov Z_{d+1, d+\tilde t} \ov Z_{d + \tilde t, i} = -\frac{1}{2}Z_{d+i 1} \ov Z_{\tilde t 1} \ov Z_{i+d, \tilde t}$ ($\tilde t = 1, \dots, d$): these simplify with the trinomials $\frac{1}{2}Z_{s1} \ov Z_{st} \ov Z_{t1}$ ($s= d+1, \dots, 2d, t=1, \dots, d$) obtained above, and then we see that also trinomials of type (I) disappear.

As for trinomials of type (III), that is $-Z_{si} \ov Z_{sj} \ov Z_{ji}$, we have that since $i,j \leq d$ then it must be $i=1$, that is $-Z_{s1} \ov Z_{sj} \ov Z_{j1}$: now, if $s \leq d$ then $j=1$ and $Z_{j1} = Z_{11} = 0$, so we are left with $-Z_{s1} \ov Z_{sj} \ov Z_{j1}$, $j \leq d, s = d+1, \dots, 2d$.

Now, we have that $\Delta_1$ contains only the trinomials of the kind $\frac{1}{2} Z_{t1} \ov Z_{ts} \ov Z_{s1}$ for $t,s = 1, \dots, 2d$: if the index $t \leq d$ then, by the same argument as above, we must have $s =1$ (otherwise $Z_{ts} = 0$), but then $Z_{s1} = Z_{11} = 0$; if $t, s$ are both greater than $d$, then by the symmetries of $Z$, we must have $t=d+1$ (otherwise $Z_{ts}=0$) and then $Z_{t1} = Z_{d+1,1} = 0$. So we are left with the trinomials $\frac{1}{2} Z_{t1} \ov Z_{ts} \ov Z_{s1}$ with $t= d+1, \dots, 2d$ and $s= 1, \dots, d$. These are exactly, up to the factor $\frac{1}{2}$, the opposite of the trinomials left in $\Delta_d$, so we see that these trinomials cancel in the expansion of the potential if and only if $c_1 = 2 c_d$, where $c_1$ and $c_d$ denote respectively the coefficient in front of $\Delta_1$ and $\Delta_d$. 

Now we prove that all the remaining cases with at least two black nodes are non Bochner. This will end the proof of Theorem \ref{thm main}.

\subsubsection{$G = SO(2d), Sp(d), SO(2d+1)$, $\Pi_M = \{e_k-e_{k+1}, e_r-e_{r+1} \} \ (k < r)$}
In all the cases but $G = SO(2d), r = d-1$, we see that the trinomial 
$$Z_{d+r,k} \ov Z_{d+r,r+1} \ov Z_{r+1,k} = z_{-e_r-e_k} \ov z_{-e_{r+1}-e_r} \ov z_{e_{r+1}-e_k}$$ 
appears in $\Delta_k$ and not in $\Delta_r$.
Indeed, it appears only as a trinomial of type (I) (with coefficient $\frac{1}{2}$) both in $\Delta_k$ and in $\Delta_r$, but in $\Delta_r$ it simplifies with the other trinomial of type (1) given by $\frac{1}{2} Z_{d+k,r} \ov Z_{d+k,d+r+1} \ov Z_{d+r+1,r}$ (which does not appear in $\Delta_k$) which equals $$\frac{1}{2}(-z_{-e_r-e_k}) (-\ov z_{e_{r+1}-e_k}) (-\ov z_{-e_{r+1}-e_r})$$ in $SO(2d)$, $SO(2d+1)$ and 
$$\frac{1}{2}(z_{-e_r-e_k}) (-\ov z_{-e_{r+1}-e_k}) (\ov z_{-e_{r+1}-e_r})$$ in $Sp(d)$. 

The same argument applies in the case $G = SO(2d), r = d-1$ to the trinomial $Z_{2d-1,1} \ov Z_{2d-1,2d} \ov Z_{2d,1} = - z_{-e_1-e_{d-1}} \ov z_{e_d-e_{d-1}} \ov z_{-e_{1}-e_d}$.

\subsubsection{$\Sp(d)$, $\Pi_M \supseteq \{e_r-e_{r+1}, 2 e_{d}\}$}\label{subsub Sp}

Assume that the last two admissible minors are $\Delta_r$ and $\Delta_d$: we are going to show that the trinomial 
$$Z_{d+1,1} \ov Z_{d+1,d} \ov Z_{d1} = z_{-2e_1} \ov z_{-e_1-e_d} \ov z_{e_d-e_1}$$ 
appears in the expansion of $\Delta_d$ (where it does not simplify with other trinomials) while it vanishes in the other admissible minors.
Indeed, on the one hand this trinomial appears as a trinomial of type (I) in the statement of Lemma \ref{lem trinomial} in all the admissible minors, with $i=1$, $s=d+1$, $t=d$; moreover, in all these minors we have also the trinomial of type (I) given by the choice $i=1$, $s=d+1$, $t= 2d$, that is $Z_{d+1,1} \ov Z_{d+1, 2d} \ov Z_{2d,1} = z_{-2e_1} (-\ov z_{e_d-e_1}) \ov z_{-e_d-e_1}$ which simplifies with $Z_{d+1,1} \ov Z_{d+1,d} \ov Z_{d1} $.
On the other hand, the trinomial $Z_{d+1,1} \ov Z_{d+1,d} \ov Z_{d1}$ appears only in the admissible minor $\Delta_d$ as a trinomial of type (III) in the statement of Lemma \ref{lem trinomial}, and as such it does not simplify with any other trinomial.

\subsubsection{$SO(2d)$, $\Pi_M \supseteq \{e_r-e_{r+1}, e_{d-1}+e_{d} \}$ and $SO(2d+1)$, $\Pi_M \supseteq \{e_r-e_{r+1}, e_{d} \}$, \ $r>1$}\label{subsub SO}

Assume that the last two admissible minors are $\Delta_r$ and $\Delta_d$: we are going to show that the trinomial 
$$Z_{d+1,2} \ov Z_{d+1,d} \ov Z_{d2} = z_{-e_1-e_2} \ov z_{-e_1-e_d} \ov z_{e_d-e_2}$$ 
appears in the expansion of $\Delta_d$ (where it does not simplify with other trinomials) while it vanishes in the other admissible minors.
Indeed, on the one hand this trinomial appears as a trinomial of type (I) in the statement of Lemma \ref{lem trinomial} in all the admissible minors, with $i=2$, $s=d+1$, $t=d$; moreover, in all these minors we have also the trinomial of type (I) given by the choice $i=1$, $s=d+2$, $t= 2d$, that is $Z_{d+2,1} \ov Z_{d+2, 2d} \ov Z_{2d,1} = (-z_{-e_1-e_2}) (-\ov z_{e_d-e_2}) (-\ov z_{-e_d-e_1})$, which simplifies with $Z_{d+1,2} \ov Z_{d+1,d} \ov Z_{d2} $.
On the other hand, the trinomial $Z_{d+1,2} \ov Z_{d+1,d} \ov Z_{d2}$ appears only in the admissible minor $\Delta_d$ as a trinomial of type (III), that is $-Z_{si} \ov Z_{sj} \ov Z_{ji}$ with $s=d+1, i=2, j=d$ and as such it does not simplify with any other trinomial.

\subsubsection{$G = SO(2d+1), \Pi_M = \{ e_1 - e_{2}, e_{d} \}$ } 
In this case we see that the trinomial 
$$Z_{2d,1} \ov Z_{2d,2d+1} \ov Z_{2d+1,1} = (-z_{-e_1-e_d}) (-\ov z_{-e_{d}}) \ov z_{-e_{1}}$$ 
appears in $\Delta_1$ and not in $\Delta_d$.
Indeed, it appears only as a trinomial of type (I) (with coefficient $\frac{1}{2}$) both in $\Delta_1$ and in $\Delta_d$, but in $\Delta_d$ it simplifies with the other trinomial of type (I) given by $\frac{1}{2} Z_{d+1,d} \ov Z_{d+1, 2d+1} \ov Z_{2d+1,d} = \frac{1}{2} z_{-e_1-e_d} (-\ov z_{-e_{d}}) \ov z_{-e_{1}}$ (which does not appear in $\Delta_1$)

\subsubsection{Any $G$, $\Pi_M$ contains at least three nodes of type $e_j-e_{j+1}$} 
More precisely, let us assume that $e_j-e_{j+1}, e_q-e_{q+1}, e_r-e_{r+1}$ are the first three black nodes, so that $\Delta_j, \Delta_q, \Delta_r$ are the first three admissible minors. Let $c_j, c_q, c_r$ the coefficients in front of these trinomials in the expansion of the potential.
We are going to show that there is no choice of the coefficients $c_j, c_q, c_r > 0$ such that the trinomials
\begin{equation}\label{trinomialAnyG1}
Z_{rj} \ov Z_{rq} \ov Z_{qj} = z_{e_r - e_j} \bar z_{e_r - e_q} \bar z_{e_q - e_j}
\end{equation}
and
\begin{equation}\label{trinomialAnyG2}
Z_{r+1,j} \ov Z_{r+1,q} \ov Z_{qj} = z_{e_{r+1} - e_j} \bar z_{e_{r+1} - e_q} \bar z_{e_q - e_j}
\end{equation}
both disappear in the expansion of (\ref{diastasisSENZAlog}).
Indeed, in $\Delta_j$ the trinomial (\ref{trinomialAnyG1}) appears only as a trinomial of type (I), so with coefficient $\frac{1}{2}$; in $\Delta_q$ it appears only as a trinomial of type (I) and (III), so with coefficient $\frac{1}{2} - 1 = -\frac{1}{2}$; and in $\Delta_r$ and all the following admissible minors it appears in all types (I)-(IV), so with coefficient $0$. Summing up, we have that (\ref{trinomialAnyG1})  appears in the expansion of (\ref{diastasisSENZAlog}) with coefficient $\frac{c_j}{2} - \frac{c_q}{2}$.
Similarly, in $\Delta_j$ the trinomial (\ref{trinomialAnyG2}) appears only as a trinomial of type (I), so with coefficient $\frac{1}{2}$; and both in $\Delta_q$ and in $\Delta_r$ it appears as a trinomial of type (I) and (III), so with coefficient $\frac{1}{2} - 1 = -\frac{1}{2}$ in bot minors;  and all the following admissible minors it appears in all types (I)-(IV), so with coefficient $0$. Summing up, we have that (\ref{trinomialAnyG2})  appears in the expansion of the potential with coefficient $\frac{c_j}{2} - \frac{c_q}{2} - \frac{c_r}{2}$.
Then, in order for both (\ref{trinomialAnyG1}) and (\ref{trinomialAnyG2}) to disappear we must have $\frac{c_j}{2} - \frac{c_q}{2}=0 $ and $\frac{c_j}{2} - \frac{c_q}{2} - \frac{c_r}{2}$, which is impossible since this implies $c_r = 0$.
The proof of Theorem \ref{thm main} is done.


\begin{thebibliography}{ABBR}
\bibitem{Ak} E. Akyildiz, {\em On the factorization of the Poincar\'e polynomial: a survey}, Serdica Math. J. 30 (2004), 159-176.
\bibitem{A-P} D. V. Alekseevsky, A. M. Perelomov, {\em Invariant Kaehler-Einstein metrics on compact homogeneous
spaces}, Funct. Anal. Appl., 20 (3) (1986) 171-182.
\bibitem{note} C. Arezzo, A. Loi,
{A note on K\"ahler-Einstein metrics and Bochner's coordinates}, {\em Abh. Math. Sem. Univ. Hamburg} 74 (2004), 49--55.
\bibitem{arv} A. Arvanitoyeorgos, {\em Geometry of flag manifolds}, International Journal of Geometric Methods in Modern Physics Vol.3, Nos. 5, 6 (2006), 957-974.
\bibitem{B} A. L. Besse, {\em Einstein manifolds}, Springer-Verag (1987).
\bibitem{BFR} M. Bordemann, M. Forger, H. Roemer, {\em Homogeneous Kaehler manifolds: paving the way towards new supersymmetric sigma models}, Commun. Math. Phys. 102, 605-647 (1986).
\bibitem{boc} S. Bochner, {\em Curvature in Hermitian metric}, Bull. Amer. Math. Soc. 53 (1947), 179-195. 
\bibitem{ca} E. Calabi, {\em Isometric Imbeddings of Complex Manifolds}, Ann. of Math. 58 (1953), 1-23. 
\bibitem{helgason} S. Helgason, {\em Differential geometry, Lie groups and symmetric spaces}, Pure and Applied Mathematics, vol. 80, Academic Press Inc. [Harcourt Brace Jovanovich Publishers], New York, 1978.
\bibitem{knapp} A. W. Knapp, {\em Lie groups beyond an introduction}, Progress in Mathematics, 140, Birkhaeuser 1996.
\bibitem{diastherm} A. Loi, {\em Calabi's diastasis function for Hermitian symmetric spaces},
Diff. Geom. Appl. 24 (2006),  311-319.
\bibitem{diastexp} A. Loi, R. Mossa, {\em The diastatic exponential of a symmetric space}, Math. Z. 268 (2011), 3-4, 1057-1068.
\bibitem{berezin} A. Loi, R. Mossa, {\em Berezin quantization of homogeneous bounded domains}, Geom. Dedicata 161 (2012), 119-128.
\bibitem{remrkhom} A. Loi, R. Mossa, {\em Some remarks on homogeneous \K\  manifolds, Geometriae dedicata. 179 (2015), no. 1, 377-383}.
\bibitem{LMZ} A. Loi, R. Mossa, F. Zuddas {\em The Log-term of the disk bundle over a homogeneous Hodge manifold}, Ann. Global Anal. Geom. (2016), in press. 
\bibitem{diastbal} R. Mossa, {\em A note on diastatic entropy and balanced metrics}, J. Geom. Phys. 86 (2014), 492-496.
\bibitem{diastent} R. Mossa,  {\em Diastatic entropy and rigidity of complex hyperbolic manifolds}, Complex Manifolds 3 (2016), 186-192.
\bibitem{ruan} W. D. Ruan, \emph{Canonical coordinates and Bergmann metrics}, Comm. in Anal. 
and Geom. (1998), 589--631.
\bibitem{tas} H. Tasaki,
{\em The cut locus and
the diastasis of a Hermitian
symmetric space of compact type},
Osaka J. Math. 22 (1985), 863-870.
\bibitem{tian4} G. Tian, \emph{On a set of polarized K\"ahler metrics on algebraic manifolds}, J. Diff. Geometry 32 (1990), 99-130.
\end{thebibliography}
\end{document}